\documentclass[letterpaper, 11pt]{amsart}
\pdfoutput=1

\usepackage{amsmath}
\usepackage{graphics}
\usepackage[english]{babel}
\usepackage{ifthen}
\usepackage{ifpdf}
\usepackage{enumerate,amsfonts,amscd,amsthm,amssymb,cite} 
\usepackage{graphicx}
\usepackage{url}

\usepackage[pdftex, plainpages = false, 
		 pdfstartview={FitH},                
		 pdfpagelayout = SinglePage,
		 bookmarks = true, bookmarksopen = true,
                 bookmarksnumbered = true,
                 linktocpage, 
                 colorlinks = true, linkcolor = blue,
                 urlcolor  = blue, citecolor = red,
                 anchorcolor = green, hyperindex = true,
                 hyperfigures]{hyperref} 
\title[Generalization of Zaslavsky's Theorem]{On a Generalization of Zaslavsky's Theorem for Hyperplane Arrangements}
\author[P\, Deshpande ]{Priyavrat Deshpande}
\address{Chennai Mathematical Institute \\ Chennai\\ India}
\email{pdeshpande@cmi.ac.in} 


\def\ds{\displaystyle}
\def\ts{\textstyle}

\newcommand{\nn}{\nonumber}
\newcommand{\ol}{\overline}
\newcommand{\Poin}{\operatorname{Poin}}
\newcommand{\rank}{\operatorname{rank}}
\numberwithin{equation}{section}

\theoremstyle{plain}
\newtheorem{theorem}{Theorem}[section]
\newtheorem{lem}[theorem]{Lemma}
\newtheorem{cor}[theorem]{Corollary}
\newtheorem{prop}[theorem]{Proposition}

\theoremstyle{definition}
\newtheorem{defn}[theorem]{Definition}
\newtheorem{ex}[theorem]{Example}
\newtheorem{xca}[theorem]{Exercise}
\newtheorem{conj}[theorem]{Conjecture}

\theoremstyle{remark}
\newtheorem{rem}[theorem]{Remark}

\newcommand{\bt}[1]{\begin{theorem}\label{#1}}
\newcommand{\bc}[1]{\begin{cor}\label{#1}}
\newcommand{\bl}[1]{\begin{lem}\label{#1}}
\newcommand{\bp}[1]{\begin{prop}\label{#1}}
\newcommand{\be}[1]{\begin{ex}\label{#1}}
\newcommand{\bd}[1]{\begin{defn}\label{#1}}
\newcommand{\brem}[1]{\begin{rem}\label{#1}}
\newcommand{\bx}[1]{\begin{xca}\label{#1}}
\newcommand{\bcon}[1]{\begin{conj}\label{#1}}

\newcommand{\et}{\end{theorem}}
\newcommand{\ec}{\end{cor}}
\newcommand{\el}{\end{lem}}
\newcommand{\ep}{\end{prop}}
\newcommand{\ee}{\end{ex}}
\newcommand{\ed}{\end{defn}}
\newcommand{\exc}{\end{xca}}
\newcommand{\erem}{\end{rem}}
\newcommand{\econ}{\end{conj}}

\newcommand{\bpr}{\begin{proof}}
\newcommand{\epr}{\end{proof}}

\newcommand{\inte}{\operatorname{Int}}
\def\A  {\mathcal{A}}
\def\FA {\mathcal{F(A)}}
\def\Ch {\mathcal{C(\A)}}

\def \F {\mathcal{F}}

\def \P {\mathcal{P}}
\def \D {\mathcal{D}}

\def \s {\mathcal{S}}
\def \L {\mathcal{L}}

\def\R{\mathbb{R}}

\def\C{\mathbb{C}}
\def\Z{\mathbb{Z}}

\begin{document}

%title page information
%\title[Generalization of Zaslavsky's Theorem]{On a Generalization of Zaslavsky's Theorem for Hyperplane Arrangements}
%\author{Priyavrat Deshpande}
%\address{Chennai Mathematical Institute\\ SIPCOT IT Park, Siruseri\\ Chennai, India}
%\email{pdeshpande@cmi.ac.in}

\begin{abstract}
We define arrangements of codimension-$1$ submanifolds in a smooth manifold which generalize arrangements of hyperplanes. When these submanifolds are removed the manifold breaks up into regions, each of which is homeomorphic to an open disc. The aim of this paper is to derive formulas that count the number of regions formed by such an arrangement. We achieve this aim by generalizing Zaslavsky's theorem to this setting. We show that this number is determined by the combinatorics of the intersections of these submanifolds. \end{abstract}

\thanks{This work is a part of author's doctoral thesis. The author would like to thank his supervisor Graham Denham for his support. Sincere thanks to Thomas Zaslavsky for discussions in Toronto, warm hospitality in Binghamton and especially for pointing out his paper on topological dissections. The paper has benefited from the inputs of the anonymous referee. The author wishes to thank the referee for a careful reading and helpful suggestions.}
\maketitle

% Introduction

Consider the problem of counting the number of pieces into which a topological space is divided when finitely many of its subspaces are removed. In the literature this is referred to as the topological dissection problem. This problem has a long history in combinatorial geometry. In 1826, Steiner considered the problem of counting the pieces of a plane cut by a finite collection of lines, circles etc. In 1901, Schl\"afli obtained a formula for counting the number of regions in a Euclidean space when it is cut by hyperplanes in general position. Subsequently many mathematicians studied various aspects and generalizations of this problem. We refer to \cite[Chapter 18]{grunbaum67} and \cite{gr72} for more information on related problems.\par 

An arrangement of real hyperplanes is a finite collection of hyperplanes in a finite dimensional real affine space. The complement of the union of these hyperplanes is disconnected. An arrangement stratifies the ambient space into open polyhedra called \textit{faces}. The top-dimensional faces are called \textit{chambers}. Zaslavsky discovered a counting formula for the number of chambers which depends on the intersection data of the arrangement in \cite[Theorem A]{zas75}. To be precise, he showed that the number of chambers is equal to the magnitude of the characteristic polynomial of the arrangement evaluated at $-1$. The characteristic polynomial is a Tutte-Grothendieck invariant of the associated intersection lattice \cite[Section 4A]{zas75}. He also showed that the $f$-polynomial (the generating polynomial of the face-counting numbers) is an evaluation of the M\"obius polynomial (which is related to the Tutte polynomial) of the intersection lattice. \par

A generalization of hyperplane arrangements, called arrangements of submanifolds, was introduced in \cite{deshpande_thesis11}. Such an arrangement is a finite collection of locally flat, codimension-$1$ submanifolds of a given manifold. The complement of the union of these submanifolds is disconnected. The aim of this paper is to generalize Zaslavsky's result to submanifold arrangements. We prove that the number of connected components of the complement of a submanifold arrangement is determined by the intersections of these submanifolds. Our result is motivated by the techniques used in \cite{ehr09}. The authors generalize Zaslavsky's formula to toric arrangements.\par 

We should point out that our main result (Theorem \ref{rthm2}) also follows from \cite[Corollary 2.2]{zas77}. However, the results are obtained independently and the techniques used here are different.\par 

The paper is organized as follows. We start Section \ref{sec:prelim} with a quick review of some combinatorial notions. In Section \ref{subsec1p1} we gather relevant background material on the theory of valuations on a poset and the Euler characteristic. In Section \ref{hyparr} we recall Zaslavsky's theorem for hyperplane arrangements. The new object of study, an arrangement of submanifolds, is introduced in Section \ref{submarr}. In Section \ref{sec:NumberOfChambers} we first define a generalization of the characteristic polynomial. We then establish a formula that combines the geometry and combinatorics of the intersections and counts the number of chambers. We also compare Zaslavsky's proof in \cite{zas77} with ours. In Section \ref{sec:TheFVector} we look at some particular cases of manifolds and derive formulas for the number of faces. Finally, in Section \ref{end} we raise some questions for further research.

% Section 1

\section{Preliminaries} \label{sec:prelim}
In this section we list the combinatorial notions we need and fix the notation. We assume the reader's familiarity with the basic concepts related to partially ordered sets and lattices; \cite[Chapter 3]{stan97} is the main reference.\par 

Let $\P = (\P, \leq)$ denote a poset. Unless stated otherwise we will always assume that $\P$ is finite. A subset $C$ of $\P$ is called a \textit{chain} if any two elements of $C$ are comparable. The \textit{length} $l(C)$ of a finite chain is defined as $|C| - 1$. A poset $\P$ is called \textit{graded} if every maximal chain of $\P$ has the same length. For such posets there is a unique \textit{rank function} $r\colon \P\to \{ 0, 1, \dots, n\}$ such that $r(x) = 0$ if $x$ is a minimal element and $r(y) = r(x) + 1$ if $y$ covers $x$. The \textit{rank} of a graded poset $\P$, denoted by  $\rank(\P)$, is defined as $\hbox{max} \{l(C)~|~ C \hbox{~is a chain of~} \P \}$. In this paper we assume that all posets are graded of some finite rank and contain a unique minimal element $\hat{0}$. \par

A notion that encodes combinatorial as well as topological information about a poset $\P$ is the M\"obius function $\mu\colon \P\times \P\to \Z$. It is defined as follows: 
\begin{align*}
\mu(x, y) = \begin{cases}
0 &\mbox{ if } y < x,\\ 
1 &\mbox{ if } y = x, \\ 
\ds -\sum_{x\leq z< y} \mu(x, z) &\mbox{ if } x< y. 
\end{cases}
\end{align*}

The M\"obius function of a poset is used to obtain inversion formulas. 
\bl{lem1ch1sec1} Let $\P$ be a poset and let $f, g\colon \P\to \C$. Then 
\[g(y) = \sum_{x\leq y} f(x) \]
if and only if \[f(y) = \sum_{x\leq y} \mu(x, y) g(x).\] \el

\bd{defchar}
For a poset $\P$ of rank $n$ the \emph{characteristic polynomial} of $\P$ is defined as:
\begin{align*}
 p(t) &:= \ds\sum_{x\in \P} \mu(\hat{0}, x) t^{n - r(x)}.
\end{align*} \ed

\subsection{Valuations on lattices}\label{subsec1p1}
We now review valuations on lattices and explain how the Euler characteristic is defined combinatorially. The main references for this material are \cite{rota71} and \cite[Chapter 2]{kla97}. \par

Let $\D$ be a family of subsets of a finite set $S$ such that $\D$ is closed under finite unions and finite intersections. Such a family is a distributive lattice in which the partial ordering is given by the inclusion of subsets, the empty set is the least element, $S$ is the greatest element while the meet and join are defined by intersection and union of subsets, respectively. 

\bd{rdef1}An $\R$-\emph{valuation} on $\D$ is a function $\nu\colon \D\to \R$, satisfying 
\begin{align}
	\ts\nu(A\cup B) &= \nu(A) + \nu(B) - \nu(A\cap B), \label{req1}\\
	\nu(\emptyset) &= 0.\label{req2}
\end{align}
By iterating the identity (\ref{req1}) we get the \emph{inclusion-exclusion principle} for $\nu$, namely
\begin{equation}\begin{split}
	\ts\nu(A_1\cup\cdots\cup A_n) = &\sum_i \nu(A_i) - \sum_{i<j}\nu(A_i\cap A_j) \\ & \quad + \sum_{i<j<k} \nu(A_i\cap A_j\cap A_k)-\cdots. \label{req3}\end{split}
\end{equation} \ed

\begin{theorem}[Rota \cite{rota71}]\label{c7s1t1}
A valuation on a finite distributive lattice $\D$ is uniquely determined by the values it takes on the set of join-irreducible elements of $\D$, and these values can be arbitrarily assigned.
\end{theorem}

With this theorem we are now in a position to define the Euler characteristic. 

\bd{c7s1d1}The \emph{Euler characteristic} of a finite distributive lattice $\D$ is the unique valuation $\chi$ such that $\chi(x) = 1$ for all join-irreducible elements $x$ and $\chi(\hat{0}) = 0$. \ed

Just like a measure, a valuation on a lattice can be used to construct abstract integrals. We explain it below. This will be used to show that the characteristic polynomial of an arrangement is an integral with respect to the Euler characteristic.\par    
For a subset $A$ of $S$, the \textit{indicator function} of $A$, denoted by $I_A$, is the function on $S$ defined by setting $I_{A}(x) = 1$ if $x\in A$ and $0$ otherwise.  A $\D$-\emph{simple} function $f\colon S\to \R$ is the following finite linear combination: 
\begin{equation}
 f = \ds\sum_{i=1}^k r_i I_{A_i}  \label{req4}
\end{equation}
where $r_i\in \R$ and $A_i\in \D$. 
The set of all $\D$-simple functions forms a ring under pointwise addition and multiplication. A subset $\L$ of $\D$ is called a \emph{generating set} if it is closed under finite intersections and if every element of $\D$ can be expressed as a finite union of members of $\L$. Using the inclusion-exclusion formula it can be shown that every $\D$-simple function can be rewritten as a linear combination 
\begin{equation}
f = \sum_{i=1}^m s_i I_{B_i} \label{req5}
\end{equation}
where each $B_i\in \L$. An $\R$-valued function $\nu$ on $L$ is called a \emph{valuation} on $\L$ provided that $\nu$ satisfies the identities (\ref{req1}) and (\ref{req2}) for all sets $A, B\in \L$ whenever $A\cup B\in \L$ (recall that $\L$ is closed under intersections). Note that, since $\L$ need not be closed under unions, identities (\ref{req1}) and (\ref{req3}) do not make sense in general. Given a valuation $\nu$ on $\L$ one could try and extend it to the whole of $\D$. Such an extension need not be well defined since an element of $\D$ can be expressed as a union of elements of $\L$ in more than way. \par 

Given a valuation $\nu$ on $\L$ and a $\D$-simple function $f$ define the \emph{integral} of $f$ with respect to $\nu$ as
\begin{equation}
	\int f d\nu = \sum_{i=1}^m s_i \nu(B_i). \label{req6}
\end{equation}
The extension of $\nu$ from $\L$ to $\D$ and the existence of the integral with respect to $\nu$ are equivalent (Groemer's integral theorem \cite[Theorem 2.2.7]{kla97}).

\subsection{Hyperplane arrangements} \label{hyparr}
Hyperplane arrangements arise naturally in geometric, algebraic and combinatorial settings, such as finite dimensional projective spaces or vector spaces defined over a field of any characteristic. Here we formally define hyperplane arrangements and the combinatorial data associated with them in a setting that is most relevant to our work. 

\bd{def1} A (real) \emph{arrangement of hyperplanes} is a finite collection $\A = \{H_1,\dots,H_k\}$ of affine hyperplanes in $\R^l$, $l\geq 1$. \ed

The \textit{rank of an arrangement} is the largest dimension of the subspace spanned by the normals to hyperplanes in $\A$. We call $\A$ a \emph{central} arrangement if all the hyperplanes pass through a common point; otherwise we call $\A$ a \emph{non-central} arrangement. An \textit{essential arrangement} is one whose rank is $l$. For a subset $F$ of $\R^l$, the \emph{restriction} of $\A$ to $F$ is the subarrangement $\A_F := \{H\in\A ~|~ F\subseteq H\}$. The hyperplanes of $\A$ induce a stratification of $\R^l$; the components of each stratum are open polyhedra and are called the \emph{faces} of $\A$. \par

There are two posets associated with $\A$, namely, the face poset and the intersection poset, which contain important combinatorial information about the arrangement. 

\bd{def2} The \emph{intersection poset} $L(\A)$ of a hyperplane arrangement $\A$ is defined as the set of all nonempty intersections of hyperplanes, including $\R^l$ itself as the empty intersection, ordered by reverse inclusion.
\ed

The rank of each element is the codimension of the corresponding intersection. In general $L(\A)$ is a (meet) semilattice; it is a lattice if and only if the arrangement is central.  

\bd{def22}The \emph{characteristic polynomial} of an arrangement $\A$ is:
\[p(\A,t) := \sum_{X\in L}\mu(\R^l, X)\cdot t^{\dim(X)}. \] \ed

\bd{def3} The \emph{face poset} $\FA$ of $\A$ is the set of all faces ordered by topological inclusion: $F\leq G$ if and only if $F\subseteq\overline{G}$. \ed

The set of all chambers (top-dimensional faces) is denoted by $\Ch$. As the complement of the hyperplanes in $\R^l$ is disconnected, a natural question is whether the number of chambers depends on the intersection data. Zaslavsky in his fundamental treatise \cite{zas75} studied the relationship between the intersection lattice of an arrangement and the number of chambers. He developed the enumeration theory for hyperplane arrangements by exploiting the combinatorial structure of the intersection lattice. His main result is as follows: 

\begin{theorem}[Theorem A \cite{zas75}]\label{them1sec1} Let $\A$ be a hyperplane arrangement in $\R^l$ with $L(\A)$ as its intersection poset and let $p(\A, t)$ be the associated characteristic polynomial. Then the number of its chambers is 
\[ |\Ch| = \sum_{X\in L(\A)} |\mu(\R^l, X)| = (-1)^l p(\A, -1).\] \et

% Section 2

\section{Submanifold arrangements}\label{submarr}

In this section we propose a generalization of arrangements of hyperplanes. In order to achieve this generalization we isolate the following characteristics of a hyperplane arrangement: 
\begin{enumerate}
\item[(1)] there are finitely many codimension-$1$ subspaces, each of which separates $\R^l$ into two components, 
\item[(2)] there is a stratification of $\R^l$ into open polyhedra, 
\item[(3)] the face poset of this stratification has the homotopy type of $\R^l$. 
\end{enumerate}

Any reasonable generalization of hyperplane arrangements should possess these properties. Since smooth manifolds are locally Euclidean they are obvious candidates for the ambient space. In this setting we can study arrangements of codimension-$1$ submanifolds that satisfy certain nice conditions; for example, locally, we would like our submanifolds to behave like hyperplanes. 

\subsection{Locally flat submanifolds}
We start by generalizing property (1) mentioned above. Throughout this paper a manifold always has empty boundary. Our focus is on the codimension-$1$ smooth submanifolds that are \emph{embedded as a closed subset} of a finite dimensional smooth manifold. This type of submanifold behaves much like a hyperplane. The following are some of its well known (separation) properties. 

\bl{lem1news1c3}
If $X$ is a connected $l$-manifold and $N$ is a connected $(l-1)$-manifold embedded in $X$ as a closed subset, then $X\setminus N$ has either $1$ or $2$ components. If in addition $H_1(X, \Z_2) = 0$ then $X\setminus N$ has two components. \el

\begin{proof} The lemma follows from the following exact sequence of pairs in mod $2$ homology:
\[ H_1(X, \Z_2)\to H_1(X, X\setminus N, \Z_2)\to \tilde{H}_0(X\setminus N, \Z_2) \to \tilde{H}_0(X, \Z_2) = 0.\qedhere  \]\end{proof} 

\bd{def1news1c3} A connected codimension-$1$ submanifold $N$ in $X$ is said to be \emph{two-sided}
if $N$ has a neighborhood $U_N$ such that $U_N\setminus N$ has two connected components; otherwise $N$ is said to be \emph{one-sided}. A submanifold $N$ is \textit{locally two-sided} in $X$ if each $x\in N$ has an arbitrarily small connected neighbourhood $U_x$ such that $U_x\setminus (U_x\cap N)$ has two components. In general a disconnected codimension-$1$ submanifold is (locally) two-sided if each of its connected components is (locally) two-sided. Moreover a submanifold \emph{separates} $X$ if its complement has $2$ components. \ed

Note that being two-sided is a local condition. For example, a point in $S^1$ does not separate $S^1$, however, it is two-sided. The following corollaries follow from the definitions and the Lemma \ref{lem1news1c3}. 

\bc{cor1s1c3}
Every codimension-$1$ submanifold $N$ is locally two-sided in $X$. \ec

\bc{cor2news1c3} If $X$ is an $l$-manifold and $N$ is an $(l-1)$ manifold embedded in $X$ as a closed subset, and if $H_1(N, \Z_2) \cong 0$, then $N$ is two-sided. \ec

An $n$-manifold $N$ contained in an $l$-manifold $X$ is \emph{locally flat} at $x\in N$ if there exists a neighborhood $U_x$ of $x$ in $X$ such that $(U_x, U_x\cap N)\cong (\R^l, \R^n).$ An embedding $f\colon N\to X$ such that $f(N)\subseteq X$ is said to be \emph{locally flat} at a point $x\in N$ if $f(N)$ is locally flat at $f(x)$. Embeddings and submanifolds are \emph{locally flat} if they are locally flat at every point.\par 
It is necessary to consider the locally flat class of submanifolds. Otherwise one could run into pathological situations. For example, the Alexander horned sphere is a (non-flat) embedding of $S^2$ inside $S^3$ such that the connected components of its complement are not even simply connected (see \cite[Page 65]{rushing73}).\par 
Locally flat submanifolds need not intersect like hyperplanes. A simple example comes from the non-Pappus arrangement of $9$ pseudolines in $\R^2$. Corresponding to this arrangement there is a rank $3$ oriented matroid which realizes an arrangement of pseudocircles in $S^2\subseteq \R^3$. Now consider cones over these pseudocircles. The cones are pseudoplanes in $\R^3$. Each pseudoplane is homeomorphic to a $2$-dimensional subspace. However there is no homeomorphism of $\R^3$ to itself mapping them onto a hyperplane arrangement. For more on pseudo-arrangements see \cite[Chapter 5]{ombook99}. \par 
We introduce a notion that will guarantee hyperplane-like intersections of submanifolds. But first some notation.
Let $\A= \{N_1, \dots, N_k \}$ be a collection of locally flat, codimension-$1$ submanifolds of $X$. For every $x\in X$ and an open neighbourhood $V_x$ (homeomorphic to $\R^l$) of $x$ let $\A_x := \{N\cap V_x ~|~ x\in N\in \A \}$. By $\bigcup \A_x$ we mean the union of elements of $\A_x$.

\bd{ch3def0} Let $X$ be a manifold of dimension $l$. Let $\A = \{N_1,\dots,N_k\}$ be a collection of codimension $1$, locally flat submanifolds of $X$. We say that these submanifolds have a \emph{locally flat intersection} 
if for every $x\in X$ there exists an open neighbourhood $V_x$ and a homeomorphism $\phi\colon V_x\to \R^l$ such that $(V_x, \bigcup \A_x)\cong (\R^l, \bigcup \A')$ where $\A'$ is a central hyperplane arrangement in $\R^l$ with $\phi(x)$ as a common point. \ed

\subsection{Cellular stratifications}
Now we generalize properties (2) and (3). Let $\A = \{N_1,\dots, N_k\}$ be a finite collection of codimension-$1$ submanifolds of $X$ having locally flat intersections. Let $\L$ denote the set of all non-empty intersections and $\L^d$ be the subset of codimension-$d$ intersections. We have $\bigcup\L^0 = X$ and $\bigcup\L^1 = \bigcup_{i=1}^k N_i$. For each $d\geq 0$ Consider the following subsets of $X$.  
\begin{align*}
	\s^d(X) &= \bigcup \L^d\setminus \bigcup \L^{d+1}.
\end{align*} 
Note that each $\s^i(X)$ may be disconnected and that $X$ can be expressed as the disjoint union of these connected components. We want these components to define a `nice' stratification of $X$. We introduce the language of \textit{cellular stratified spaces} developed in \cite{tamaki01} in order to achieve the generalization. Recall that a subset $A$ of a topological space $X$ is \emph{locally closed} if every point $x\in A$ has a neighbourhood $U$ in $X$ with $A\cap U$ closed in $U$. 

\bd{defn3s1}Let $X$ be a topological space and $\P$ be a poset. A \textit{stratification} of $X$ indexed by $\P$ is a surjective map $\sigma\colon X\to \P$ satisfying the following properties:
\begin{enumerate}
	\item For $p\in \P$, $e_p := \sigma^{-1}(p)$ is connected and locally closed.
	\item For $p, q\in \P,~ e_p\subseteq \ol{e_q} \iff p\leq q$.
	\item $e_p \cap \ol{e_q} \neq \emptyset \Longrightarrow e_p \subseteq \ol{e_q}$.
\end{enumerate}
The subspace $e_p$ is called the \textit{stratum with index $p$}.\ed

One can verify that the boundary of each stratum, $\partial e_p = \ol{e_p} - e_p$, is itself a union of strata. Such a stratification gives a decomposition of $X$, i.e., it is the disjoint union of strata. 
The indexing poset $\P$ will be called the \emph{face poset}. Let $A$ be a subspace of $X$. If the restriction $\sigma|_{A}$ is a stratification $(A, \sigma|_{A})$ then it is called a \textit{stratified subspace} of $(X, \sigma)$. If, in addition, $A$ is a union of cells of the stratification then 
it is called a \textit{strict stratified subspace}. \par 

It is now easy to check that the components of $\s^i(X)$ define a stratification of $X$ when ordered by inclusion. However it is not desirable to consider arbitrary stratifications. For example, consider two non-intersecting longitudinal circles in the $2$-torus $S^1\times S^1$; there are two codimension-$0$ strata and two codimension-$1$ strata. The resulting face poset does not have the homotopy type of the torus. We need to focus on stratifications such that the strata are cells. We make this precise.

\bd{defn4s1}
A \textit{globular $n$-cell} is a subset $D$ of $D^n$ (the unit $n$-disk) containing $\mathrm{Int} (D^n)$. We call $D\cap \partial(D^n)$
the boundary of $D$ and denote it by $\partial D$. The number $n$ is called the \textit{globular dimension} of $D$.\ed
\bd{defn5s1} Let $X$ be a Hausdorff space. A \textit{cellular stratification} of $X$ is a pair $(\sigma, \Phi)$ of a stratification $\sigma$ and a collection of continuous maps, $\Phi = \{\phi_p : D_p\to \ol{e_p}~|~ p\in \P\}$, called \textit{cell structures}, satisfying the following conditions:
\begin{enumerate}
	\item Each $D_p$ is a globular cell and $\phi_p|_{\mathrm{Int}(D^n)}\colon \mathrm{\inte}(D^n)\to e_p$ is a homeomorphism.
	\item $\phi_p(D_p) = \ol{e_p}$ and $\phi_p$ is a quotient map for every $p\in \P$. 
	\item For each $n$-cell $e_p$ the boundary $\partial e_p$ contains cells of dimension $i$ for every $0\leq i\leq n-1$.
\end{enumerate} A \textit{cellularly stratified space} is a triple $(X, \sigma, \Phi)$ where $(\sigma, \Phi)$ is a cellular stratification. \ed
We assume that all our cellularly stratified spaces (CS-spaces for short) have a finite number of cells. Consequently the cell structure is CW (i.e., each cell meets only finitely many other cells and X has the weak topology determined by the union of cells). 
A cell $e_p$ is said to be \emph{regular} if the structure map $\phi_p$ is a homeomorphism. A cell complex is regular if all its cells are regular.

\bd{defn6s1} Let $X$ be a cellularly stratified space. $X$ is called \textit{totally normal} if for each $n$-cell $e_p$,
\begin{enumerate}
	\item there exists a structure of regular cell complex on $S^{n-1}$ containing $\partial D_p$ as a strict cellular stratified subspace of $S^{n-1}$ and
	\item for any cell $D$ in $\partial D_p$, there exists a cell $e_q\in \partial e_p$ such that $D$ is homeomorphic to $D_q$ and the characteristic map $\phi_q = \phi_p|_{\ol{D}}$.
\end{enumerate} \ed

In particular, the closure of each $k$-cell contains at least one cell of dimension $i$, for every $i\leq k-1$. We now give an example of a cellularly stratified space that is not totally normal. Consider $\mathrm{Int}(D^2)\cup \{(1, 0)\}$; its boundary does not contain a $1$-cell.

\bd{defn7s1} The \emph{face category} of a cellularly stratified space $(X, \sigma, \Phi)$ is denoted by $\F(X)$. The objects of this category are the cells of $X$. For each pair $p \leq q$, define $\F(X)(e_p, e_q)$ (i.e., the set of all morphisms) to be the set of all maps $b: D_p\to D_q$ such that $\phi_q\circ b = \iota \circ \phi_p$. \ed

Recall that a small category $\F$ is said to be \textit{acyclic} if, for any pair of distinct objects $x, y\in\mathrm{Ob}(\F)$, either $\F(x, y)$ or $\F(y, x)$ is empty, and for any object $x$, $\F(x, x)$ consists of the identity morphism. \par 

A poset $(\P, \leq)$ can be regarded as a an acyclic category as follows. The objects of this category are the elements of $\P$. There is a unique morphism from $x$ to $y$ if $x\leq y$. In fact, one can think of acyclic categories as a  generalization of posets. For any acyclic category $\F$ there exists a unique partial order $\preceq$ on its objects given by:
\[\F(x, y)\neq \emptyset \Rightarrow x\preceq y. \]
$(\mathrm{Ob}(\F), \preceq)$ is called the \textit{underlying poset} of $\F$. 
Let $\F$ be an acyclic category such that there is at most one morphism between any two objects. In that case $\F$ and its underlying poset are equivalent as acyclic categories. This is precisely what it will mean when we say that an acyclic category is a poset. \par 

\begin{lem}[Lemma 4.2\cite{tamaki01}]\label{lem0s1}
If $X$ is a cellularly stratified space then its face category $\F(X)$ is acyclic. When $X$ is regular, $\F(X)$ is equivalent to the face poset of $X$. \el

For CS-spaces the underlying poset coincides with the (classical) face poset. The reason we are using the language of CS-spaces is the following. 

\bl{lemn1s1} Let $\A$ be an arrangement of hyperplanes in $\R^l$. Then the stratification induced by the hyperplanes in $\A$ defines a structure of  totally normal, regular cellularly stratified space on $\R^l$. \el

%\bpr It is enough to characterize the cell structure maps. Each stratum of the arrangement is a relatively open convex polyhedron. The intersection of the closure of such a stratum with the open unit disc is a sector which is (PL) homeomorphic to a globular cell of an appropriate dimension. Checking the other conditions of the definition is now straightforward.\epr

As a quick example consider the arrangement of coordinate axes in $\R^2$. For the closed face $e = \{(x,y) ~|~ x, y\geq 0 \}$ the subspace $D = \mathrm{Int}(D^2)\cup \{(x,y)\in S^1~|~ x < 0\}$ can serve as the globular cell. \par 

Now we state a result that will explain the need for total normality. Recall that the face poset of a regular CW complex has the homotopy type of the complex. This property is also shared by the totally normal CS-spaces.  

\begin{lem}[Corollary  4.17\cite{tamaki01}] \label{lemn2s1} For a totally normal cellularly stratified space $X$, the nerve of its face category embeds in $X$ as a strong deformation retract. \el

The closure of each cell in a totally normal CS-space is homeomorphic to a polyhedral complex \cite[Section 3.3]{tamaki01}. In a nutshell, the strata and the face category of a totally normal CS-space behave similarly to that of hyperplane arrangements. If we assume that such a structure is induced due to a finite collection of submanifolds then the properties (2) and (3) mentioned at the beginning of this section are generalized.

\subsection{Definition and examples}

The desired generalization of hyperplane arrangements is the following: 
\bd{def31}Let $X$ be a connected, smooth, real manifold of dimension $l$. An \textbf{arrangement of submanifolds} is a finite collection 
$\A = \{N_1,\dots, N_k\}$ of codimension-$1$ smooth submanifolds in $X$ such that:
	\begin{enumerate}
		\item the $N_i$'s have locally flat intersection,
		\item the stratification induced by the intersections of $N_i$'s defines the structure of a totally normal cellularly stratified space on $X$.
	\end{enumerate}
\ed

Let us look at how we can associate combinatorial data to such an arrangement and at a few examples. 

\bd{def32}The \emph{intersection poset}, denoted by $L(\A)$, is the set of connected components of nonempty intersections of $N_i$'s, ordered by reverse inclusion. The rank of each element in $L(\A)$ is defined to be the codimension of the corresponding intersection.\ed
By convention $X$ is the intersection of no submanifolds, hence it is the smallest member of $L(\A)$. Note that in general this poset need not be a lattice. The use of connected components of intersections is not a new idea, see for example \cite{zas77} and more recently (in the case of toric arrangements) \cite{ehr09, moci_tutte_2009}. 

\bd{def33}The cells of the stratification will be called the \emph{faces} of the arrangement. The \textit{face category} of the arrangement is the face category of the induced stratification and is denoted by $\FA$. Codimension-$0$ faces are called \textit{chambers} and the set of all chambers is denoted by $\Ch$. \ed

Hyperplane arrangements are obvious examples of these submanifold arrangements. Here are some examples of arrangements in spheres.

\be{ex41} Let $X$ be the circle $S^1$, a smooth $1$-manifold. The codimension-$1$ submanifolds are points in $S^1$. Consider the arrangement $\A = \{ p\}$ of one point. There are two strata, one $0$-cell $p$ and one $1$-cell $A$. In this case there are two lifts of the cell structure map for $p$. Hence in the face category there are two morphisms. The face category is not equivalent to its underlying poset. The intersection poset consists of only two elements.\par  
Figure \ref{s1pt1fig} shows this arrangement with a graph depicting the geometric realization of the face category. The two edges correspond to the two morphisms from $p$ to $A$. Note that the geometric realization has the homotopy type of $S^1$.
\begin{figure}[!ht]
\begin{center}  \includegraphics[scale=0.3,clip]{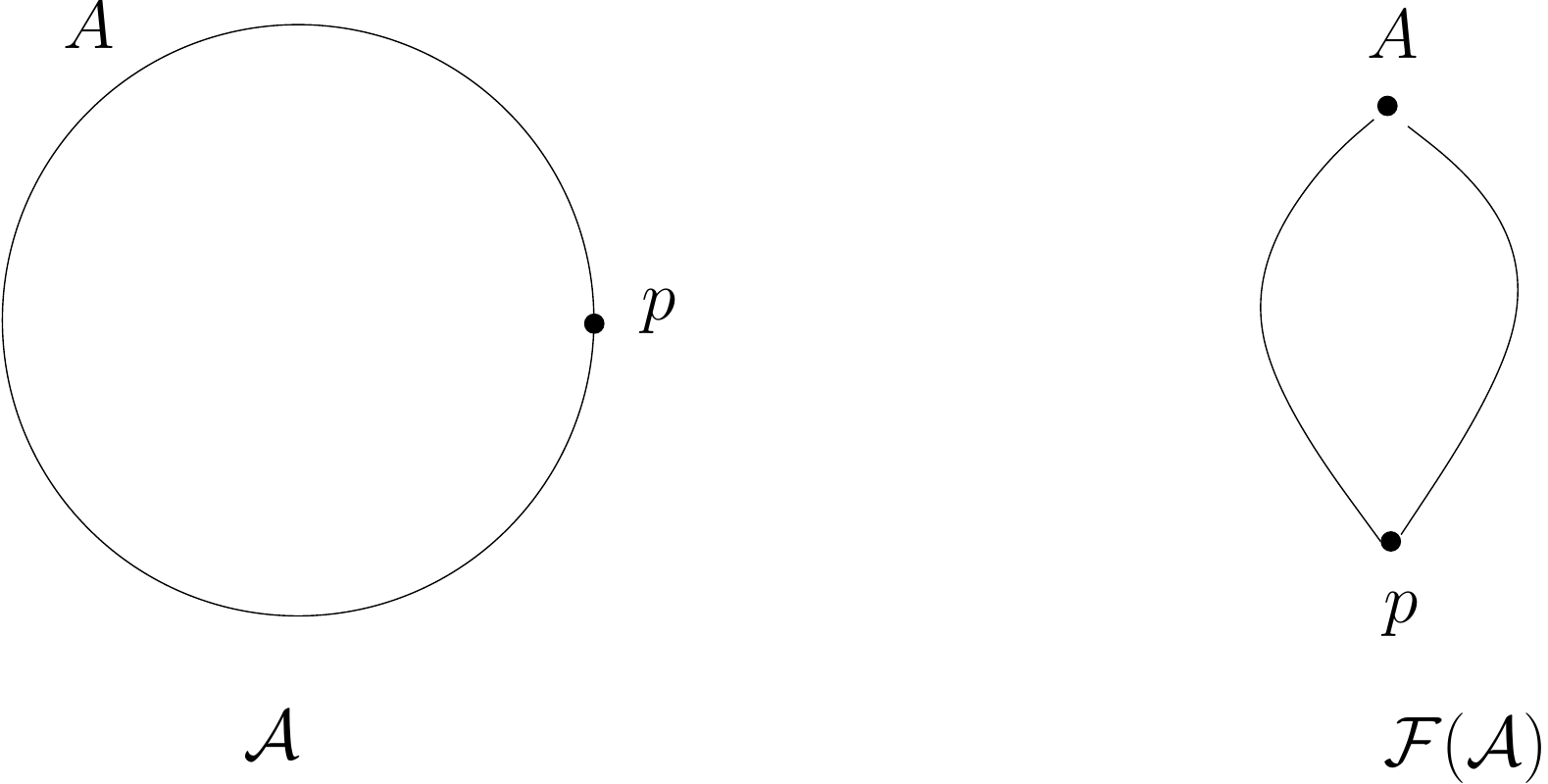} \end{center}
    \caption{Arrangement of $1$ point in a circle.}     \label{s1pt1fig}     \end{figure}
\ee

\be{ex42} Consider the arrangement $\A = \{p, q\}$ of $2$ points in $S^1$. Here the induced cellular stratification is regular, hence the face category is equivalent to the face poset. Figure \ref{figdef1} shows this arrangement and the Hasse diagrams of the face poset and the intersection poset. 
\begin{figure}[!ht]
\begin{center}  \includegraphics[scale=0.35,clip]{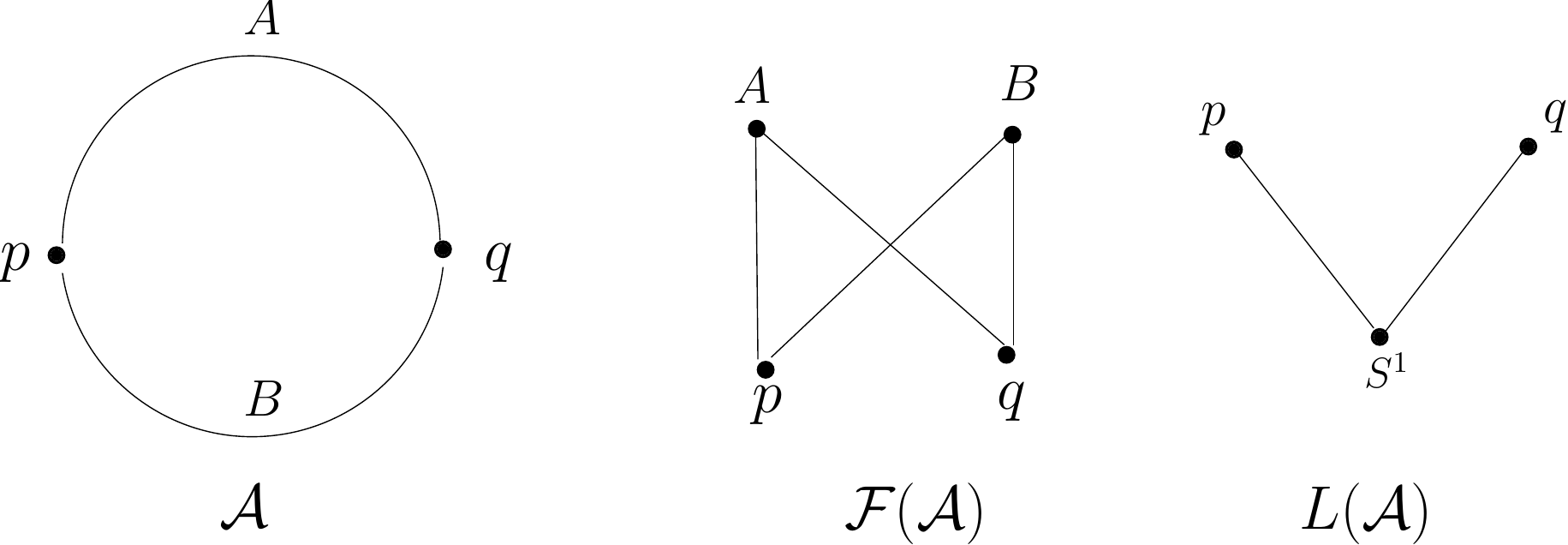} \end{center}
    \caption{Arrangement of $2$ points in a circle.}     \label{figdef1}     \end{figure}
\ee

\be{ex43}As a $2$-dimensional example consider an arrangement of $2$ great circles $N_1, N_2$ in $S^2$. The stratification consists of two $0$-cells, four $1$-cells and four $2$-cells. Figure \ref{sphere01} shows this arrangement and the Hasse diagram of the intersection poset. 
\begin{figure}[!ht]
\begin{center}\includegraphics[scale=0.35,clip]{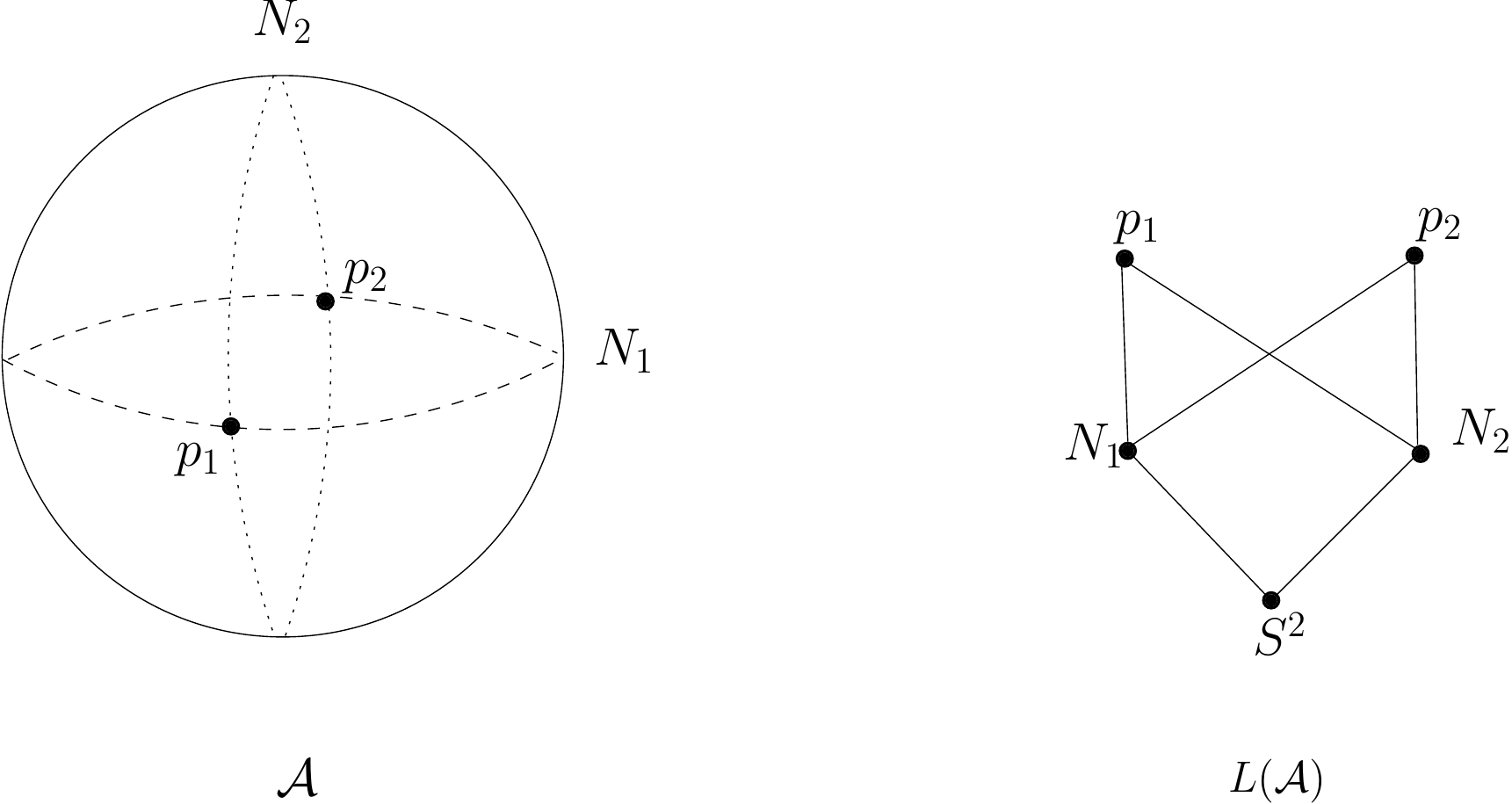}\end{center}
    \caption{Arrangement of $2$ circles in a sphere.} \label{sphere01}     \end{figure} \ee

\be{ex44}Now consider the $2$-torus $T^2 = \R^2/\Z^2$. A codimension-$1$ submanifold of $T^2$ is the image of a straight line (with rational slope) passing through the origin under the projection $\R^2\to T^2$. Let $\A$ be the arrangement in $T^2$ obtained by projecting the lines $x = -2y, y = -2x, y = x$. These toric hyperplanes intersect in $3$ points, namely $p_1 = (0,0)$, $p_2 = (1/3, 1/3)$ and $p_3 = (2/3, 2/3)$. There are $6$ chambers. Figure \ref{torus} shows the arrangement ($T^2$ is considered as the quotient of the unit square) with the associated intersection poset.
\begin{figure}[!ht]
\begin{center}\includegraphics[scale=0.35,clip]{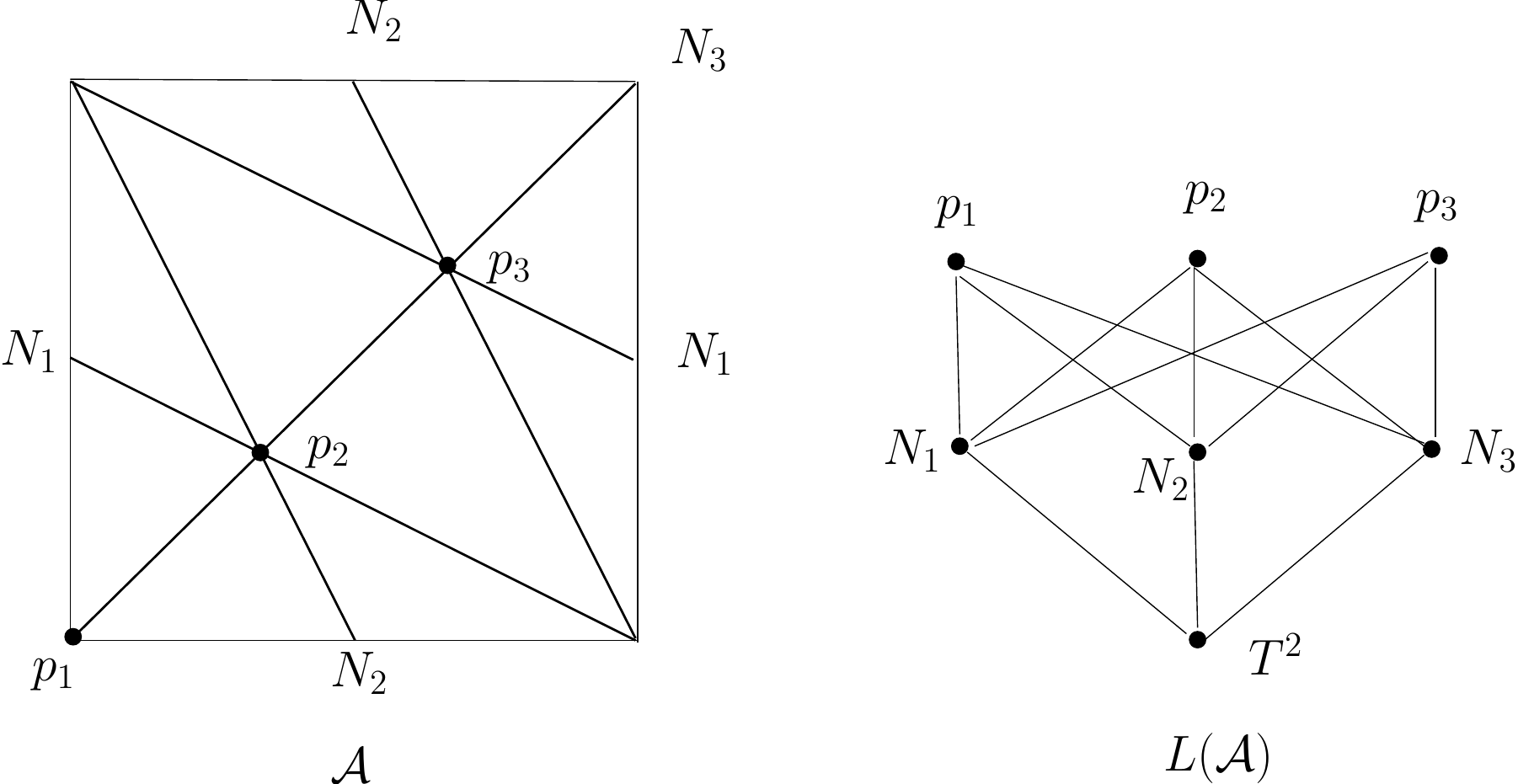}\end{center}
    \caption{Arrangement in a $2$-torus.} \label{torus}     \end{figure}
\ee

Example \ref{ex44} is an example of a toric arrangement. Recently there has been a surge of interest in toric arrangements (see for example \cite{dantonio_salvetti_2011, ehr09, jim2011870, moci_tutte_2009}).  

% Section 3

\section{The Chamber Counting Formula} \label{sec:NumberOfChambers}

Let $\A$ be an arrangement of submanifolds of a smooth $l$- manifold $X$. The problem at hand is to count the number of chambers of $\A$. Let $\D$ be the lattice of sets generated by elements of $L(\A)$ and chambers through unions and intersections. We start by generalizing the characteristic polynomial.\par

Define the \emph{Poincar\'e polynomial with compact support} of a topological space $A$ as 
\[\Poin_c(A,t) := \sum_{i\geq 0}\rank(H^i_c(A,\Z)) t^i\]
where $H^i_c$ is the cohomology with compact supports.

\bl{rlem1} The function $\nu\colon \D\to \Z[t]$ defined by 
\[ \nu(A) = \Poin_c(A,t),~\forall A\in \D \]
is a valuation on $\D$. \el

\bpr The first step is to find a generating set for $\D$. Let ${G}$ denote the set of all possible finite unions of faces of $\A$. Clearly ${G}$ is closed under intersections. Each member of $L(\A)$ is a union of finitely many faces. Since all the chambers are faces any finite union (or intersection) of members of $L(\A)$ and $\Ch$ can be expressed as a finite union of faces. This proves that ${G}$ is a generating set for $\D$. \par
Now we have to show that $\nu$ defines a valuation on $G$. This is clear because $\nu(\emptyset) = 0$ and all the non-empty faces are disjoint, open topological cells. The inclusion-exclusion identities (\ref{req3}) are satisfied by $\nu$ because of the excision property. Hence as a consequence of \emph{Groemer's integral theorem} $\nu$ extends to the whole of $\D$ (see \cite[Theorem 2.2.7]{kla97}).  \epr

For each $Y\in L ~(= L(\A))$, define \[f(Y) = Y \setminus \bigcup_{\substack{Z \subseteq Y\\ Z\in L}} Z.\]
The set $f(Y)$ is a union of faces and $\{f(Y) ~|~ Y\in L\}$ is a disjoint collection. Hence
\[ X = \coprod_{\substack{X\leq Y\\ Y\in L}}f(Y).\] 
Consequently for the indicator functions we have the following equation:
\[ I_X = \sum_{X\leq Y} I_{f(Y)}.\]
Therefore by M\"obius inversion 
\begin{equation}
I_{f(X)} = \sum_{X\leq Y}\mu(X, Y) I_Y.\label{req8}
\end{equation}
Since $f(X)$ is the union of all the chambers, $I_{f(X)}$ is a $\D$-simple function. As $\nu$ extends uniquely to a valuation on $\D$ it defines an integral on the algebra of $\D$-simple functions. Integrating $I_{\Ch}$ with respect to $\nu$ we get 
\begin{align}
\int I_{\Ch} d\nu &= \sum_{Y\in L} \mu(X,Y) \nu(Y) \nn \\
	     &= \sum_{Y\in L} \mu(X,Y) \Poin_c(Y,t).
\end{align}

\bd{c7s2d1} Let $\A$ be an arrangement of submanifolds in an $l$-manifold $X$ and $L$ be the associated intersection poset. The \emph{generalized characteristic polynomial} of $\A$ is \[\ds p(\A, t) := \sum_{Y\in L} \mu(X,Y) \Poin_c(Y,t). \] \ed

\bd{c7s2d2} The \emph{combinatorial Euler characteristic} $\kappa$ of a finite CW complex $P$ is defined as:
\[\kappa(P) = \begin{cases}\chi(\hat{P}) - 1 & \hbox{~if $P$ is not compact,}\\ \chi(P) & \hbox{~if $P$ is compact}\end{cases} \]
where $\hat{P}$ is the one-point compactification of $P$ and $\chi$ is the usual Euler characteristic. \ed

\begin{rem}\label{remce} If $F$ is a face of an arrangement $\A$ then $\kappa(F) = (-1)^{\dim F}$. \end{rem}

Note that this is not a new notion, this is just a topological description of Definition \ref{c7s1d1}. To give an intrinsic topological description of the combinatorial Euler characteristic for arbitrary spaces is not an easy job. The theory of o-minimal structures has to be used in order to define valuations and integrals on arbitrary spaces, which is beyond the scope of this paper. However, the above notion is a topological invariant and it satisfies the Euler relation, that is, the number of even-dimensional cells minus the number of odd-dimensional cells is equal to the value of $\kappa$.\par 
The following lemma is now clear.

\bl{c7s2l1} $\kappa$ defines an $\R$-valuation on $\D$. \el
We can now generalize Zaslavsky's theorem. The following theorem also follows from \cite[Corollary 2.2]{zas77}. We discuss this point towards the end of this section.

\bt{rthm2}
Let $\A$ be an arrangement of submanifolds in an $l$-manifold $X$. Then the number of chambers is given by 
\begin{equation} |\Ch| = (-1)^l\sum_{Y\in L} \mu(X,Y)\kappa(Y). \label{req10} \end{equation} \et

\bpr First note that  $\kappa$ and $\nu|_{t = -1}$ agree on every element of $G$. Consequently, they also agree on every member of $\D$. Hence, we have  
\begin{equation} \kappa(\Ch) = \int I_{\Ch} d\kappa = \int I_{\Ch} d\nu |_{t = -1}.\label{req9}\end{equation}
The set $\Ch$ is a disjoint union of chambers, each of which is homeomorphic to an open ball of dimension $l$. From Remark \ref{remce} and the above equations we get
\begin{align} 
|\Ch| &= (-1)^l \kappa(\Ch)\notag \\
&= (-1)^l \int I_{\Ch} d\nu |_{t = -1}\notag \\ 
					 &= (-1)^l\sum_{Y\in L} \mu(X,Y) \Poin_c(Y,-1)\notag \\
					 &= (-1)^l\sum_{Y\in L} \mu(X,Y)\kappa(Y). \qedhere
\end{align} 
\end{proof}
Hence the number of chambers of an arrangement not only depends on the combinatorics of the intersections but also on their geometry.

\begin{cor}[Theorem A \cite{zas75}] \label{rcor1} Let $\A$ be a hyperplane arrangement in $\R^n$. The number of chambers of this arrangement is equal to $(-1)^n \chi(\A, -1)$.  \ec
\bpr Every member of the intersection poset in this case is homeomorphic to an open ball; hence, $\Poin_c(Y, t) = t^{\dim Y}$ for every $Y\in L(\A)$. The result follows from the observation that $\int I_{\Ch} d\nu = \chi(\A, t)$.\epr

Next we recover the result in the context of toric arrangements (see \cite[Section 3]{ehr09} for the technical definition).

\begin{cor}[Theorem 3.11 \cite{ehr09}] \label{rcor2} For a toric hyperplane arrangement $\A$ in a torus $T^n$ that subdivides the torus into  $n$-cells, the number of chambers is given by \[\ds (-1)^n \sum_{\dim Y = 0}\mu(T^n, Y).\] \ec

\bpr Note that toric hyperplanes are homeomorphic to $T^{n-1}$ and that their intersections are disjoint unions of subtori. If $Y\in L$ then $\Poin_c(Y, t) = (1+t)^k$ for $0\leq k\leq n$. Substituting it in equation (\ref{req10}) we get the formula. 
\epr

In case of Example \ref{ex43} the generalized characteristic polynomial is $p(\A, t) = t^2 - 2t + 1$ and $p(\A, -1) = 4$, which is the number of chambers. For the toric arrangement of Example \ref{ex44} we have $p(\A, t) = t^2 - t + 4$ and $p(\A, -1) = 6$, the number of chambers.\par 

Recently Forge and Zaslavsky have considered arrangements of topological hyperplanes in \cite{zas09}; these are slightly more general than arrangements of pseudohyperplanes. However, the counting formula remains the same because the connected components of the complement are cells.\par 

Now consider a more general situation where $X$ is a topological space and $\A$ is a finite collection of subspaces that are removed from $X$. Let $L$ denote the poset consisting of $X$ and the connected components of all possible finite intersections of members of $\A$ ordered by reverse inclusion. The topological dissection problem asks whether it is possible to express the number of connected components of the complement in terms of the combinatorics of the intersection poset $L$ (dissection of $X$ intuitively means that $X$ is expressed as a union of pairwise disjoint subspaces). Let $C_1,\dots,C_m$ denote the connected components of $X\setminus \bigcup\A$. In this context Theorem \ref{rthm2} takes the following form:

\bt{c7s2t3} If the combinatorial Euler characteristic $\kappa$ is a valuation on the lattice of sets generated by $L\cup \{C_1,\dots,C_m\}$ through unions and intersections then
\[\sum_{j=1}^m\kappa(C_j) = \sum_{Y\in L} \mu(X,Y)\kappa(Y). \] \et

The above statement is referred to as the \emph{fundamental theorem of dissection theory} in \cite[Theorem 1.2]{zas77}. In a nutshell, the number of connected components of the complement (i.e., the combinatorial Euler characteristic of the complement) depends on a condition on the intersections. 
%This condition turns out to be the Euler relation. 
Moreover if every face of the dissection is a finite, disjoint union of open cells then the combinatorial Euler characteristic is a valuation  \cite[Lemma 1.1]{zas77}. \par 

The approach we took to prove Theorem \ref{rthm2} can be traced back to a paper of Ehrenborg and Readdy \cite{ehr98}, where they determined the characteristic polynomial of a subspace arrangement. They explicitly used Gromer's integral theorem to prove that the characteristic polynomial is a valuation. While studying toric arrangements with M. Slone \cite{ehr09} they generalized the previous result and proved the above mentioned Corollary \ref{rcor2}. The idea of looking at the Euler characteristic as an integral of indicator functions is due to Chen \cite{chen93}, see also \cite{chen00}. \par

Finally we compare our strategy with that of Zaslavsky's in \cite{zas77}. At the very foundation of both the strategies lies the idea of using the combinatorial Euler characteristic as a valuation. In order to implement this idea we have used M\"obius inversion whereas Zaslavsky used a technical property of valuations. Instead of stating his theorem we only mention an important and relevant consequence. 

\begin{theorem}[Corollary 2.2\cite{zas77}]\label{c7s2t1} Let $\phi$ be a valuation of the finite distributive lattice $\D$, and let $\P$ be a subset of $\D$ containing $\hat{0}$ and every join-irreducible element. Then for any $t\in \P$ which is not $\hat{0}$ or a join-irreducible element of $\D$,
\[ \sum_{s\in \P; s\leq t} \mu_{\P}(s, t)\phi(s) = 0.\] \et

In order to apply this theorem to the counting problem we take $\D$ to be the lattice of sets generated by the intersection poset and all the chambers. Theorem \ref{rthm2} then follows once we use the valuation $\kappa$.  

% Section 4

\section{Faces of an Arrangement} \label{sec:TheFVector}
In this section we use formula (\ref{req10}) to find the number of faces of various dimensions of an arrangement. We start with a definition.

\bd{c7s3d1}For an arrangement of submanifolds $\A$ in an $l$-dimensional manifold $X$, define its \textit{$f$-polynomial} as 
\[\ds f_{\A}(x) := \sum_{k=0}^l f_k x^{l-k}\]
where $f_k$ denotes the number of $k$-dimensional faces of the arrangement.\ed
 
For $Y\in L(\A)$ consider the set  
\[\A^Y := \{N\cap Y~|~ N\in \A \hbox{~and~} \emptyset\neq N\cap Y\neq Y\}. \]
In general $Y$ need not be a smooth manifold and the collection $\A^Y$ need not define a submanifold arrangement. However, $\A^Y$ still defines a stratification of $Y$. We refer to codimension-$0$ strata as chambers of $\A^Y$. Note that for every $Y\in L(\A)$ every chamber of $\A^Y$ is an open ($\dim Y$)-cell (in fact they are faces of $\A$). Hence Theorem \ref{c7s2t3} can be used to count the number of chambers of $\A^Y$, which in turn gives us the numbers $f_k$. Such a dissection is referred to as a \emph{dissection with properly cellular faces} in \cite{zas77}.

\bt{rthm3} Let $X$ be a smooth, real manifold of dimension $l$ and $\A$ be an arrangement of submanifolds. Then the numbers $f_k$ are given by 
\begin{equation}
f_k = \sum_{\dim Y = k}(-1)^k\sum_{\substack{Z\in L\\ Y\leq Z}}\mu(Y,Z) \kappa(Z). \label{c7s3eq1}\end{equation} \et

\bpr As $\FA = \{\mathcal{C}(\A^Y)~ |~ Y\in L \}$, the number of $k$-faces of $\A$ is given by
\[f_k = \sum_{\substack{Y\in L(\A) \\ \dim Y = k }} |\mathcal{C}(\A^Y)|. \]
Use (\ref{req10}) to substitute for $|\mathcal{C}(\A^Y)|$ in the above formula and note that 
\[L(\A^Y) = \{Z\in L(\A)~|~ Y\leq Z \}. \]
The $f$-polynomial takes the form
\[\ds f_{\A}(x) = (-1)^l\sum_{Z\in L(\A)}\kappa(Z)\sum_{Y\leq Z}\mu(Y, Z)(-x)^{l - \dim Y}. \qedhere\] \epr 

Recall that a finite lattice $L$ is called \emph{atomic} if every $x\in L$ is a join of atoms (the elements that cover $\hat{0}$).  A finite lattice is called \emph{semimodular} if when $x, y$ both cover $x\wedge y$, then $x\vee y$ covers both $x$ and $y$. A \emph{geometric} lattice is both atomic and semimodular (see for example \cite[Chapter3]{stan97}). It is a well known fact that the intersection lattice of a central hyperplane arrangement is geometric.  The intersection poset of a submanifold arrangement need not be a lattice but every interval in it is a geometric lattice, as is proved in the following lemma.

\bl{rlem2}Let $\A$ be an arrangement of submanifolds in an $l$-manifold $X$. Then every interval of the intersection poset $L(\A)$ is a geometric lattice. \el
\bpr Consider an interval $[Y, Z]$ in $L(\A)$ such that  $\dim Y = i$ and $\dim Z = j$. There exist an open set $V$ in $X$ and a coordinate chart $\phi$ such that $\phi(V\cap Y)$ is homeomorphic to an $i$-dimensional subspace of $\R^l$. Moreover, $\{\phi(N\cap V) \mid N\in \A^Y \}$ is a central arrangement of hyperplanes in $\R^i$. For any $W\in [Y,Z]$ the subspace $\phi(W\cap V)$ is homeomorphic to a subspace of $\R^i$ that contains $\phi(Z\cap V)$. In particular the ($i-1$)-dimensional subspaces in $[Y, Z]$ map to hyperplanes in $\R^i$ that contain $\phi(Z\cap V)$. This correspondence gives us an essential central arrangement of hyperplanes in $\R^{i-j}$ when we quotient out $\phi(Z\cap V)$. This correspondence is also a poset isomorphism and hence $[Y, Z]$ is a geometric lattice. 
\epr
For geometric lattices the M\"obius function alternates in sign 
(see \cite[Proposition 3.10.1]{stan97}). Hence the expression for the $f$-polynomial simplifies as follows:
\begin{align}
f_{\A}(x) &= \sum_{Z\in L} \kappa(Z) \sum_{Y\leq Z} (-1)^l (-1)^{l-\dim Y} \mu(Y, Z) x^{l - \dim Y} \nn\\
		  &= \sum_{Z\in L} \kappa(Z) \sum_{Y\leq Z} (-1)^{\dim Z} |\mu(Y, Z)| x^{l - \dim Y}. \label{eqn001}
\end{align}

We say that the submanifolds of $X$ are in \emph{relative general position} if the intersection of any $i$ of the submanifolds, $i\geq 1$, is either empty or $(l-i)$-dimensional. An arrangement $\A$ of submanifolds is called \emph{simple} if the submanifolds are in relative general position. For simple arrangements, every interval of $L(\A)$ is a Boolean algebra. Hence the rank generating function for an interval $[Y, Z]$ is $(x+1)^{l-\dim Z}$ and $\mu(Y, Z) = (-1)^{\dim Y - \dim Z}$ (see \cite[Example 3.8.3]{stan97}). For simple arrangements the $f$-polynomial takes the following form:

\begin{equation}
f_{\A}(x) = \sum_{Z\in L} \kappa(Z) \sum_{Y\leq Z} (-1)^{\dim Z} (x + 1)^{l - \dim Y}. \label{eqn002}
\end{equation} 

Now we compute the face numbers $f_k$ for some particular submanifold arrangements. In the following examples $a_k$ denotes the number of rank $k$ elements in $L(\A)$. For similar calculations also see \cite{ehr09, zas77}.\par  
We start with a known example. Recall that hyperplanes are in \emph{(absolute) general position} if the intersection of $k$ hyperplanes, $1\leq k\leq l+1$, is $(l-k)$-dimensional. 

\be{rexc4} If $X \cong \R^l$ and $\A$ is an arrangement of hyperplanes then 
\[f_k = \sum_{\dim Y = k}~\sum_{Y\leq Z}|\mu(Y, Z)|.  \]
For the proof, for every element $Z\in L$ substitute $\kappa(Z) = (-1)^{\dim Z}$ in Equation (\ref{c7s3eq1}). If $\A$ is a simple arrangement then
$f_{\A}(x) = \sum_{j=0}^l a_j (x + 1)^{l - j}$ and consequently   
\[f_k = \sum_{j=0}^k a_j \binom{l-j}{l-k}. \] 
If the hyperplanes are in absolute general position then $a_j = \binom{n}{l-j}$. Substituting this in the above expression we get the formula for the number of faces of the arrangement.\ee

\be{rexc5}
Now consider an arrangement of codimension-$1$ subtori in an $l$-torus. In this case 
\[f_k = \sum_{\substack{\dim Y = k\\ \dim Z = 0\\ Y\leq Z}} |\mu(Y, Z)|.  \]
This follows from the fact that every intersection is a disjoint union of lower-dimensional tori. So $\kappa(Z) = 0$ when $\dim Z \geq 1$. If the subtori are in relative general position then the number $f_k$ takes the following simpler form
\[f_k = a_0\binom{l}{l-k}. \]
The $f$-polynomial in this case is $a_0(x + 1)^l$. \ee
See \cite{ehr09, jim2011870, moci_tutte_2009} for similar formulas and other types of enumeration problems in a torus. \par 

Let $S^l$ denote the unit sphere in $\R^{l+1}$. A subset $S$ of $S^l$ is called a \emph{hypersphere} if each connected component of $S^l\setminus S$ is homeomorphic to the $l$-ball. An arrangement of hyperspheres is a finite collection of hyperspheres in $S^l$ such that every non-empty intersection is a subsphere. As a prototypical example one can consider the intersection of a central hyperplane arrangement in $\R^{l+1}$ with $S^l$.

\be{rexc6} If $\A$ is an arrangement of hyperspheres in $S^l$ then 
\[f_k = 2\sum_{\dim Y = k}\sum_{\substack{Y\leq Z\\ \dim Z \geq 2,~ \mathrm{even}}}|\mu(Y,Z)| +\sum_{\dim Y = k}\sum_{\substack{Y\leq Z\\ \dim Z = 0}}|\mu(Y, Z)|.\]
The above formula can be obtained from Equation (\ref{c7s3eq1}) by substituting $\kappa(Z) = 2$ for even-dimensional intersections and $0$ for the odd-dimensional intersections. For simple arrangements, substituting appropriate $\kappa$ values in Equation (\ref{eqn002}) we get the following:
\[f_k = 2 \sum_{j = 2, \mathrm{even}}^{k}a_j\binom{l-j}{l-k} + a_0\binom{l}{k}. \] \ee
Zaslavsky has considered arrangements of `great spheres' in \cite{zas85}. Pakula has considered counting problems for pseudosphere arrangements in a topological sphere in \cite{pakula03}. These arrangements are more general than what we have considered above. For these arrangements the stratification obtained need not be cellular. However, the chambers of such an arrangement and all its subarrangements are assumed to have homology of a point. 

\be{rexc7} Let $X$ be the $l$-dimensional projective space and $\A$ be an arrangement of projective hyperplanes. Then substituting $\kappa(Z) = 1$ for even-dimensional intersections and $0$ for the odd we get 
\[f_k = \sum_{\dim Y = k}\sum_{\substack{Y\leq Z\\ \dim Z \geq 0,~\mathrm{even}}}|\mu(Y,Z)|. \]
If the arrangement is simple then the projective hyperplanes are in absolute general position (since parallelism is excluded) so $a_j = \binom{n}{l-j}$. Substituting this in Equation (\ref{eqn002}) we get a well-known formula 
\[f_k = \sum_{j = 0, \mathrm{even}}^{k}\binom{n}{l-j}\binom{l-j}{l-k}. \]
The formula first appeared in \cite[Section 4]{buck43}, where Buck studied counting problem for simple arrangements of hyperplanes (see also \cite[Section 5E]{zas75} and \cite[Section 18.1]{grunbaum67}). 
\ee

\subsection{Bayer-Sturmfels map}
We now describe a relationship between the intersection poset and the face poset which extends a result due to Bayer and Sturmfels \cite{bayers90} for hyperplane arrangements. Let $\F^*$ be the dual of the face poset. Define the map $\psi\colon\F^*\to L(\A)$ by sending each face to the smallest-dimensional intersection that contains the face. For oriented matroids this is the support map to its underlying matroid. The map $\psi$ is order and rank preserving as well as surjective. Hence we look at it as a map from the set of chains of $\F^*$ to the set of chains of $L(\A)$. 

\bt{rthm5} Let $c =\{Y_1 \leq Y_2 \leq\cdots\leq Y_k \},~ k\geq 2 $, be a chain in the intersection poset $L(\A)$ of an arrangement of submanifolds $\A$. Then the cardinality of the inverse image of the chain $c$ under the map $\psi$ is given by the following formula
\[|\psi^{-1}(c)| = \prod_{i=1}^{k-1}\sum_{Y_i\leq Z\leq Y_{i+1}}|\mu(Y_i, Z)|\cdot |\sum_{Y_k\leq Z} \mu(Y_k, Z)\kappa(Z)|. \]  \et

\bpr The arguments are similar to the proof of \cite[Theorem 3.13]{ehr09}. The number of ways of selecting a face $F_k$ such that $\psi(F_k) = Y_k$ is equal to the number of chambers of $\A^{Y_k}$. A face $F_{k-1}$ is in $\psi^{-1} (Y_{k-1})$ if it is a chamber of $\A^{Y_{k-1}}$ and whose closure contains the face $F_k$. The number of such faces is equal to the number of chambers in the central hyperplane arrangement whose intersection lattice is isomorphic to $[Y_{k-1}, Y_k]$. By repeating this process for all the subspaces up to $Y_1$ we get the desired formula.\epr

% Section 5

\section{Concluding Remarks}\label{end}
Usually cell enumeration problems are studied in the context of cell complexes homeomorphic to either a disc or a sphere. For example, there is  abundant literature concerning $f$-vectors of polytopes, Dehn-Sommerville equations, the $g$-theorem for simplicial polytopes etc., see \cite{grunbaum67}. In this paper we have considered cell complexes which are not necessarily homeomorphic to spheres. Studying the analogues of traditional enumeration problems (for polytopes) in this general context seems very interesting. In this section we outline some such questions. We direct the interested reader to \cite[Section 5]{ehr09} for more open questions regarding regular subdivisions of manifolds. \par

It is well known that the polar dual of a central hyperplane arrangement is a zonotope. For certain submanifold arrangements the zonotopes are replaced by metrical-hemisphere (MH for short) complexes. An MH-complex is a regular cell complex whose face poset behaves very much like that of a zonotope. For example, the $2$-cells are combinatorially equivalent to regular $2k$-gons, the MH-complex is dual to the cell structure induced by a submanifold arrangement etc. (see \cite[Chapter 3]{deshpande_thesis11}). It will be interesting to study enumeration problems for these complexes. \par

The face poset of a central hyperplane arrangement corresponds to a realizable oriented matroid. Under this correspondence the intersection lattice corresponds to the underlying matroid. The chambers of an arrangement correspond to topes of the oriented matroid. The definition of the characteristic polynomial (or more generally that of the Tutte polynomial) for central arrangements coincides with the usual definition of these polynomials for matroids.\par  

Zaslavsky's theorem has an analogue in oriented matroids. Las Vergnas proved a formula that counts the number of topes in oriented matroids, see \cite{lvergnas75, lvergnas80}. Our current work in progress concerns finding a generalization of oriented matroids. The combinatorial structure arising from the face poset of an MH-complex resembles that of an oriented matroid. Also, it is possible to find an analogue of matroids in the setting of submanifold arrangements (very similar to the multiplicity matroids discovered in \cite{moci_tutte_2009}). The aim here is to come up  with a Tutte polynomial which explains the deletion-restriction phenomenon in this general setting and also incorporates the generalized characteristic polynomial (Definition \ref{c7s2d1}). \par

One would also like to understand the combinatorial type of the chambers. It is not hard to see that closure of a chamber is a nice manifold with corners (i.e., every point has a neighborhood diffeomorphic to a neighborhood of a point in $[0, \infty)^r\times \R^l$). In some cases the face poset of a closed chamber is isomorphic to that of a convex polytope. Given an arrangement of submanifolds is it possible to characterize (or to count) the bounded chambers that are combinatorially simplices or cubes or any other type of polytope?\par 

\bibliographystyle{alpha} % bibliography style abbrv
\bibliography{zaslavskyref} % References file

\newcommand{\etalchar}[1]{$^{#1}$}
\begin{thebibliography}{BLVS{\etalchar{+}}99}

\bibitem[BLVS{\etalchar{+}}99]{ombook99}
Anders Bj{\"o}rner, Michel Las~Vergnas, Bernd Sturmfels, Neil White, and
  G{\"u}nter~M. Ziegler.
\newblock {\em {O}riented {M}atroids}, volume~46 of {\em Encyclopedia of
  Mathematics and its Applications}.
\newblock Cambridge University Press, Cambridge, second edition, 1999.

\bibitem[BS90]{bayers90}
Margaret Bayer and Bernd Sturmfels.
\newblock Lawrence polytopes.
\newblock {\em Canad. J. Math.}, 42(1):62--79, 1990.

\bibitem[Buc43]{buck43}
R.~C. Buck.
\newblock Partition of space.
\newblock {\em Amer. Math. Monthly}, 50:541--544, 1943.

\bibitem[Che93]{chen93}
Beifang Chen.
\newblock On the {E}uler characteristic of finite unions of convex sets.
\newblock {\em Discrete Comput. Geom.}, 10(1):79--93, 1993.

\bibitem[Che00]{chen00}
Beifang Chen.
\newblock On characteristic polynomials of subspace arrangements.
\newblock {\em J. Combin. Theory Ser. A}, 90(2):347--352, 2000.

\bibitem[dD12]{dantonio_salvetti_2011}
Giacomo d'Antonio and Emanuele Delucchi.
\newblock A {S}alvetti complex for toric arrangements and its fundamental
  group.
\newblock {\em Int. Math. Res. Notices}, 2012(15):3535--3566, 2012.

\bibitem[Des11]{deshpande_thesis11}
Priyavrat Deshpande.
\newblock {\em Arrangements of Submanifolds and the Tangent Bundle Complement}.
\newblock PhD thesis, The University of Western Ontario, 2011.
\newblock \url{http://ir.lib.uwo.ca/etd/154}.

\bibitem[ER98]{ehr98}
Richard Ehrenborg and Margaret~A. Readdy.
\newblock On valuations, the characteristic polynomial, and complex subspace
  arrangements.
\newblock {\em Adv. Math.}, 134(1):32--42, 1998.

\bibitem[ERS09]{ehr09}
Richard Ehrenborg, Margaret Readdy, and Michael Slone.
\newblock Affine and toric hyperplane arrangements.
\newblock {\em Discrete Comput. Geom.}, 41(4):481--512, 2009.

\bibitem[FZ09]{zas09}
David Forge and Thomas Zaslavsky.
\newblock On the division of space by topological hyperplanes.
\newblock {\em European J. Combin.}, 30(8):1835--1845, 2009.

\bibitem[Gr{\"u}67]{grunbaum67}
Branko Gr{\"u}nbaum.
\newblock {\em Convex Polytopes}.
\newblock Pure and Applied Mathematics, Vol. 16. Wiley Interscience, New York,
  1967.

\bibitem[Gr{\"u}72]{gr72}
Branko Gr{\"u}nbaum.
\newblock {\em Arrangements and Spreads}.
\newblock CBMS Reg. Conf. Ser. Math., No. 10. American Mathematical Society,
  Providence, R.I., 1972.

\bibitem[KR97]{kla97}
Daniel~A. Klain and {Gian-Carlo} Rota.
\newblock {\em Introduction to Geometric Probability}.
\newblock Lezioni Lincee. Cambridge University Press, Cambridge, 1997.

\bibitem[Law11]{jim2011870}
Jim Lawrence.
\newblock Enumeration in torus arrangements.
\newblock {\em European J. Combin.}, 32(6):870--881, 2011.

\bibitem[LV75]{lvergnas75}
Michel Las~Vergnas.
\newblock Matro\"\i des orientables.
\newblock {\em C. R. Acad. Sci. Paris S\'er. A-B}, 280:Ai, A61--A64, 1975.

\bibitem[LV80]{lvergnas80}
Michel Las~Vergnas.
\newblock Convexity in oriented matroids.
\newblock {\em J. Combin. Theory Ser. B}, 29(2):231--243, 1980.

\bibitem[Moc12]{moci_tutte_2009}
Luca Moci.
\newblock A {T}utte polynomial for toric arrangements.
\newblock {\em Trans. Amer. Math. Soc.}, 364(2):1067--1088, 2012.

\bibitem[Pak03]{pakula03}
Lewis Pakula.
\newblock Pseudosphere arrangements with simple complements.
\newblock {\em Rocky Mountain J. Math.}, 33(4):1465--1477, 2003.

\bibitem[Rot71]{rota71}
Gian-Carlo Rota.
\newblock On the combinatorics of the {E}uler characteristic.
\newblock In {\em Studies in {P}ure {M}athematics ({P}resented to {R}ichard
  {R}ado)}, pages 221--233. Academic Press, London, 1971.

\bibitem[Rus73]{rushing73}
T.~Benny Rushing.
\newblock {\em Topological Embeddings}.
\newblock Pure and Applied Mathematics, Vol. 52. Academic Press, New York,
  1973.

\bibitem[Sta97]{stan97}
Richard~P. Stanley.
\newblock {\em Enumerative {C}ombinatorics. {V}olume 1}.
\newblock Cambridge Studies in Advanced Mathematics, Vol. 49. Cambridge
  University Press, Cambridge, 1997.

\bibitem[Tam]{tamaki01}
Dai Tamaki.
\newblock Cellular stratified spaces {I}: Face categories and classifying
  spaces.
\newblock arXiv:1106.3772v3 [math.AT].

\bibitem[Zas]{zas75}
Thomas Zaslavsky.
\newblock Facing up to arrangements: face-count formulas for partitions of
  space by hyperplanes.
\newblock {\em Mem. Amer. Math. Soc.}
\newblock Vol. 1, issue 1, No. 154, 1975.

\bibitem[Zas77]{zas77}
Thomas Zaslavsky.
\newblock A combinatorial analysis of topological dissections.
\newblock {\em Advances in Math.}, 25(3):267--285, 1977.

\bibitem[Zas85]{zas85}
Thomas Zaslavsky.
\newblock Extremal arrangements of hyperplanes.
\newblock In {\em \emph{Discrete Geometry and Convexity ({N}ew {Y}ork, 1982)}},
  volume 440 of {\em Ann. New York Acad. Sci.}, pages 69--87. New York Acad.
  Sci., New York, 1985.

\end{thebibliography}

\end{document}